\documentclass{article}
\usepackage{amssymb}

\newtheorem{theorem}{Theorem}
\newtheorem{acknowledgement}[theorem]{Acknowledgement}

\input{tcilatex}

\begin{document}

\begin{center}
\textbf{On some convergences to the constant }$e$ \textbf{and improvements
of Carlemans' inequality}

\[
\]%
Cristinel Mortici$^{1}$ and Hu Yue$^{2}$\bigskip 

$^{1}$ Valahia University of T\^{a}rgovi\c{s}te, Bd. Unirii 18, 130082 T\^{a}%
rgovi\c{s}te, Romania, email: cristinel.mortici@hotmail.com

$^{2}$ College of Mathematics and Informatics, Henan Polytechnic University,
Jiaozuo City, Henan 454000, China, email: huu3y2@163.com

\[
\]
\end{center}

\[
\]

\begin{quote}
\textbf{ABSTRACT: }We present sharp inequalities related to the sequence $%
\left( 1+1/n\right) ^{n}$ and some applications to Kellers' limit and
Carlemans' inequality.
\end{quote}

\[
\]



\section{Introduction and Motivation}

The starting point of this paper is the following well-known double
inequality%
\begin{equation}
\frac{e}{2n+2}<e-\left( 1+\frac{1}{n}\right) ^{n}<\frac{e}{2n+1},\ \ \ n\geq
1.  \label{s}
\end{equation}%
This inequality was highly discussed and extended in the recent past, since
it was used to improve inequalities of Hardy-Carleman type. See for example 
\cite{yp}, \cite{sz}, \cite{y}, \cite{x}, \cite{yu}.

As (\ref{s}) is equivalent to%
\[
\frac{2n}{2n+1}<\frac{1}{e}\left( 1+\frac{1}{n}\right) ^{n}<\frac{2n+1}{2n+2}%
, 
\]%
we prove that the best approximation of the form%
\begin{equation}
\frac{1}{e}\left( 1+\frac{1}{n}\right) ^{n}\approx \frac{n+a}{n+b},\text{ \
\ as }n\rightarrow \infty  \label{na}
\end{equation}%
is obtained for $a=5/12$ and $b=11/12.$ Then we prove the following\bigskip

\textbf{Theorem 1. }\emph{For every real number }$x\in \lbrack 1,\infty ),$ 
\emph{the following inequalities hold:}%
\begin{eqnarray*}
&&\frac{x+\frac{5}{12}}{x+\frac{11}{12}}-\frac{5}{288x^{3}}+\frac{343}{%
8640x^{4}}-\frac{2621}{41\,472x^{5}} \\
&<&\frac{1}{e}\left( 1+\frac{1}{x}\right) ^{x} \\
&<&\frac{x+\frac{5}{12}}{x+\frac{11}{12}}-\frac{5}{288x^{3}}+\frac{343}{%
8640x^{4}}-\frac{2621}{41\,472x^{5}}-\frac{2621}{41\,472x^{5}}+\frac{300\,901%
}{3483\,648x^{6}}.
\end{eqnarray*}

As application, we give a new proof of the limit%
\begin{equation}
\lim_{n\rightarrow \infty }\left( \frac{\left( n+1\right) ^{n+1}}{n^{n}}-%
\frac{n^{n}}{\left( n-1\right) ^{n-1}}\right) =e.  \label{z}
\end{equation}%
This limit is also known as Keller's limit. See e.g. \cite{j}, where a
different proof of (\ref{z}) is presented.

Moreover, the estimates from Theorem 1 are strong enough to prove%
\begin{equation}
\lim_{n\rightarrow \infty }n^{2}\left( \left( \frac{\left( n+1\right) ^{n+1}%
}{n^{n}}-\frac{n^{n}}{\left( n-1\right) ^{n-1}}\right) -e\right) =\frac{e}{24%
},  \label{v}
\end{equation}%
which is a new result, according to the best of our knowledge.

Finally, improvements of Carlemans' inequality are given.

\section{The Proofs}

In order to find the best approximation (\ref{na}), we associate the
relative error sequence $w_{n}$ by the relations%
\[
\frac{1}{e}\left( 1+\frac{1}{n}\right) ^{n}=\frac{n+a}{n+b}\cdot \exp
w_{n},\ \ \ n\geq 1 
\]%
and we consider an approximation (\ref{na}) to be better when $w_{n}$
converges faster to zero. We have%
\[
w_{n}=n\ln \left( 1+\frac{1}{n}\right) -1-\ln \frac{n+a}{n+b}, 
\]%
but using a mathematical software such as Maple, we get%
\[
w_{n}=\left( -a+b-\frac{1}{2}\right) \frac{1}{n}+\left( \frac{1}{2}a^{2}-%
\frac{1}{2}b^{2}+\frac{1}{3}\right) \frac{1}{n^{2}}+\left( \frac{1}{3}b^{3}-%
\frac{1}{3}a^{3}-\frac{1}{4}\right) \frac{1}{n^{3}}+O\left( \frac{1}{n^{4}}%
\right) . 
\]%
This form can be also obtained by direct computation.

Evidently, the fastest sequence $w_{n}$ is obtained when the first two
coefficients in this expansion vanish, that is $a=5/12$ and $b=11/12.$ Our
first aim is now attained.\bigskip

\emph{Proof of Theorem 1. }The requested inequalities can be written as $f>0$
and $g<0,$ where%
\[
f\left( x\right) =x\ln \left( 1+\frac{1}{x}\right) -1-\ln \left( \frac{x+%
\frac{5}{12}}{x+\frac{11}{12}}-\frac{5}{288x^{3}}+\frac{343}{8640x^{4}}-%
\frac{2621}{41\,472x^{5}}\right) 
\]%
\[
g\left( x\right) =x\ln \left( 1+\frac{1}{x}\right) -1-\ln \left( \frac{x+%
\frac{5}{12}}{x+\frac{11}{12}}-\frac{5}{288x^{3}}+\frac{343}{8640x^{4}}-%
\frac{2621}{41\,472x^{5}}+\frac{300\,901}{3483\,648x^{6}}\right) . 
\]%
We have%
\[
f^{\prime \prime }\left( x\right) =\frac{A\left( x-1\right) }{x^{2}\left(
x+1\right) ^{2}\left( 12x+11\right) ^{2}P^{2}\left( x\right) }>0 
\]%
and%
\[
g^{\prime \prime }\left( x\right) =-\frac{B\left( x-1\right) }{x^{2}\left(
x+1\right) ^{2}\left( 12x+11\right) ^{2}Q^{2}\left( x\right) }<0, 
\]%
where%
\[
P\left( x\right)
=59\,184x^{2}-66\,708x-43\,200x^{3}+1036\,800x^{5}+2488\,320x^{6}-144\,155 
\]%
\begin{eqnarray*}
Q\left( x\right) &=&5945\,040x-5603\,472x^{2}+4971\,456x^{3}-3628\,800x^{4}
\\
&&+87\,091\,200x^{6}+209\,018\,880x^{7}+16\,549\,555
\end{eqnarray*}

\begin{eqnarray*}
A\left( x\right) &=&387\,888\,768\,\allowbreak
643\,091\,163x+1374\,068\,561\,\allowbreak 183\,884\,363\allowbreak x^{2} \\
&&+2856\,411\,438\,\allowbreak
418\,498\,368x^{3}+3861\,333\,058\,\allowbreak 156\,847\,712\allowbreak x^{4}
\\
&&+3547\,125\,026\,\allowbreak
642\,062\,080x^{5}+2242\,448\,726\,\allowbreak 942\,859\,264\allowbreak x^{6}
\\
&&+963\,345\,615\,\allowbreak 805\,707\,264x^{7}+269\,162\,452\,\allowbreak
894\,408\,704\allowbreak x^{8} \\
&&+44\,174\,729\,\allowbreak 709\,158\,400x^{9}+3234\,548\,\allowbreak
057\,702\,400\allowbreak x^{10} \\
&&+48\,685\,659\,\allowbreak 681\,079\,707
\end{eqnarray*}

\begin{eqnarray*}
B\left( x\right) &=&5495\,\allowbreak 336\,279\,092\,\allowbreak
271\,136\,793x+22\,015\,\allowbreak 820\,845\,590\,\allowbreak
210\,733\,374\allowbreak x^{2} \\
&&+52\,587\,\allowbreak 526\,363\,654\,\allowbreak
958\,754\,048x^{3}+83\,107\,\allowbreak 983\,906\,845\,\allowbreak
638\,539\,984\allowbreak x^{4} \\
&&+91\,197\,\allowbreak 790\,053\,279\,\allowbreak
643\,410\,048x^{5}+70\,886\,\allowbreak 916\,929\,730\,\allowbreak
329\,339\,904\allowbreak x^{6} \\
&&+39\,022\,\allowbreak 307\,420\,738\,\allowbreak
572\,320\,768x^{7}+14\,907\,\allowbreak 444\,982\,230\,\allowbreak
536\,515\,584\allowbreak x^{8} \\
&&+3763\,\allowbreak 807\,019\,677\,\allowbreak
591\,584\,768x^{9}+565\,\allowbreak 244\,311\,814\,\allowbreak
774\,194\,176\allowbreak x^{10} \\
&&+38\,\allowbreak 255\,330\,631\,\allowbreak
116\,390\,400x^{11}+621\,\allowbreak 810\,333\,384\,\allowbreak
191\,039\,953\allowbreak .
\end{eqnarray*}%
Evidently, $g$ is strictly concave, $f$ is strictly convex, with $f\left(
\infty \right) =g\left( \infty \right) =0,$ so $g<0$ and $f>0$ on $[1,\infty
).$ The proof is completed.$\square $

\section{Kellers' limit}

Let us rewrite Theorem 1 in the form%
\[
u\left( n\right) <\frac{1}{e}\left( 1+\frac{1}{n}\right) ^{n}<v\left(
n\right) , 
\]%
where%
\[
u\left( x\right) =\frac{x+\frac{5}{12}}{x+\frac{11}{12}}-\frac{5}{288x^{3}}+%
\frac{343}{8640x^{4}}-\frac{2621}{41\,472x^{5}} 
\]%
and%
\[
v\left( x\right) =\frac{x+\frac{5}{12}}{x+\frac{11}{12}}-\frac{5}{288x^{3}}+%
\frac{343}{8640x^{4}}-\frac{2621}{41\,472x^{5}}-\frac{2621}{41\,472x^{5}}+%
\frac{300\,901}{3483\,648x^{6}}. 
\]

We prove (\ref{z}) using policemen lemma. As the sequence%
\[
x_{n}=\frac{1}{e}\left( \frac{\left( n+1\right) ^{n+1}}{n^{n}}-\frac{n^{n}}{%
\left( n-1\right) ^{n-1}}\right) 
\]%
can be written as%
\[
x_{n}=\left( n+1\right) \frac{1}{e}\left( 1+\frac{1}{n}\right) ^{n}-n\frac{1%
}{e}\left( 1+\frac{1}{n-1}\right) ^{n-1}, 
\]%
we use Theorem 1 to obtain%
\begin{equation}
\left( n+1\right) u\left( n\right) -nv\left( n-1\right) <x_{n}<\left(
n+1\right) v\left( n\right) -nu\left( n-1\right) .  \label{t}
\end{equation}%
The extreme-side sequences are rational functions of $n$ and they tends
together to $1,$ as $n$ approaches infinity. Indeed,%
\[
\left( n+1\right) u\left( n\right) -nv\left( n-1\right) =\frac{%
2508\,226\,560n^{13}-12\,\allowbreak 959\,170\,560n^{12}+\cdots }{%
17\,418\,240n^{5}\left( n-1\right) ^{6}\left( 12n-1\right) \left(
12n+11\right) } 
\]%
and%
\[
\left( n+1\right) v\left( n\right) -nu\left( n-1\right) =\frac{%
2508\,226\,560n^{13}-10\,\allowbreak 450\,944\,000n^{12}+\cdots }{%
17\,418\,240n^{6}\left( n-1\right) ^{5}\left( 12n-1\right) \left(
12n+11\right) }. 
\]%
It results that $x_{n}$ tends to $1,$ as $n$ approaches infinity, so (\ref{z}%
) is proved.

Further, by (\ref{t}), we get%
\begin{eqnarray*}
&&n^{2}\left( \left( \left( n+1\right) u\left( n\right) -nv\left( n-1\right)
\right) -1\right)  \\
&<&n^{2}\left( x_{n}-1\right)  \\
&<&n^{2}\left( \left( \left( n+1\right) v\left( n\right) -nu\left(
n-1\right) \right) -1\right) 
\end{eqnarray*}%
and again the extreme-side sequences are rational functions of $n$ and they
tends together to $1/24,$ as $n$ approaches infinity. Indeed,%
\[
n^{2}\left( \left( \left( n+1\right) u\left( n\right) -nv\left( n-1\right)
\right) -1\right) =\frac{104\,509\,440n^{11}-539\,965\,440n^{10}+\cdots }{%
17\,418\,240n^{3}\left( n-1\right) ^{6}\left( 12n-1\right) \left(
12n+11\right) }
\]%
\[
n^{2}\left( \left( \left( n+1\right) v\left( n\right) -nu\left( n-1\right)
\right) -1\right) =\frac{104\,509\,440n^{11}+-435\,456\,000n^{10}\cdots }{%
17\,418\,240n^{4}\left( n-1\right) ^{5}\left( 12n-1\right) \left(
12n+11\right) }.
\]%
In consequence, $n^{2}\left( x_{n}-1\right) $ tends to $1/24,$ which is (\ref%
{v}).

\section{Improvements of Carlemans' inequality}

While Swedish mathematician Torsten Carleman was studying quasi-analytical
functions, he discovered an important inequality, now known as Carlemans'
inequality. If $\sum a_{n}$ is a convergent series of nonneagtive reals, then%
\begin{equation}
\sum_{n=1}^{\infty }\left( a_{1}a_{2}\cdots a_{n}\right) ^{1/n}\leq
e\sum_{n=1}^{\infty }a_{n}.  \label{e}
\end{equation}%
This inequality was proven to be of great independent interest, since many
authors improved it in the recent past.

The main tool for studying and improving (\ref{e}) was the proof of P\'{o}%
lya (see \cite{p1}-\cite{p2}), who started from AM-GM inequality in the form%
\begin{equation}
\left( a_{1}a_{2}\cdots a_{n}\right) ^{1/n}\leq \frac{c_{1}a_{1}+c_{2}a_{2}+%
\cdots +c_{n}a_{n}}{n\left( c_{1}c_{2}\cdots c_{n}\right) ^{1/n}},  \label{c}
\end{equation}%
where $c_{1},$ $c_{2},$ $...,$ $c_{n}>0.$ The proof of the following result
is based on P\'{o}lya's\ idea.\bigskip 

\textbf{Theorem 2. }\emph{Let }$a_{n}>0$ \emph{such that }$\sum a_{n}$ \emph{%
is convergent and }$c_{n}>0$ \emph{such that}%
\[
\sum_{k=1}^{\infty }\frac{1}{k\left( c_{1}c_{2}\cdots c_{k}\right) ^{1/k}}%
=l\in 
\mathbb{R}
.
\]%
\emph{Denote}%
\[
x_{n}=\sum_{k=n}^{\infty }\frac{1}{k\left( c_{1}c_{2}\cdots c_{k}\right)
^{1/k}}.
\]%
\emph{Then}%
\begin{equation}
\sum_{n=1}^{\infty }\left( a_{1}...a_{n}\right) ^{1/n}\leq
\sum_{n=1}^{\infty }c_{n}x_{n}a_{n}.  \label{cn}
\end{equation}

\emph{Proof. }Using (\ref{c}), we have%
\begin{eqnarray*}
\sum_{n=1}^{\infty }\left( a_{1}a_{2}\cdots a_{n}\right) ^{1/n} &\leq
&\sum_{n=1}^{\infty }\left( \frac{1}{n\left( c_{1}c_{2}\cdots c_{n}\right)
^{1/n}}\sum_{m=1}^{n}c_{m}a_{m}\right)  \\
&=&\sum_{n=1}^{\infty }\sum_{m=1}^{n}\frac{c_{m}a_{m}}{n\left(
c_{1}c_{2}\cdots c_{n}\right) ^{1/n}} \\
&=&\sum_{m=1}^{\infty }\sum_{n=m}^{\infty }\frac{c_{m}a_{m}}{n\left(
c_{1}c_{2}\cdots c_{n}\right) ^{1/n}} \\
&=&\sum_{m=1}^{\infty }c_{m}a_{m}\sum_{n=m}^{\infty }\frac{1}{n\left(
c_{1}c_{2}\cdots c_{n}\right) ^{1/n}} \\
&=&\sum_{m=1}^{\infty }c_{m}a_{m}\sum_{n=m}^{\infty }\left(
x_{n}-x_{n+1}\right)  \\
&=&\sum_{m=1}^{\infty }c_{m}x_{m}a_{m}.\square 
\end{eqnarray*}%
P\'{o}lya took $c_{n}=\left( n+1\right) ^{n}/n^{n-1}$ and (\ref{cn}) becomes%
\[
\sum_{n=1}^{\infty }\left( a_{1}...a_{n}\right) ^{1/n}\leq
\sum_{n=1}^{\infty }\left( 1+\frac{1}{n}\right) ^{n}a_{n}.
\]%
Now (\ref{c}) follows from $\left( 1+1/n\right) ^{n}<e.$

Almost all improvements stated in the recent past use upper bounds for $%
\left( 1+1/n\right) ^{n},$ stronger than $\left( 1+1/n\right) ^{n}<e.$ See 
\cite{yp}, \cite{j}-\cite{yu}.

We use Theorem 1 to establish the following improvement of Carlemans'
inequality.\bigskip 

\textbf{Theorem 3. }\emph{Let }$a_{n}>0$ \emph{such that }$\sum a_{n}<\infty
.$ \emph{Then}%
\[
\sum_{n=1}^{\infty }\left( a_{1}...a_{n}\right) ^{1/n}\leq
e\sum_{n=1}^{\infty }\frac{12n+5}{12n+11}a_{n}.
\]%
\emph{It also holds good}%
\[
\sum_{n=1}^{\infty }\left( a_{1}...a_{n}\right) ^{1/n}\leq
e\sum_{n=1}^{\infty }\left( \frac{12n+5}{12n+11}-\varepsilon _{n}\right)
a_{n},
\]%
\emph{where}%
\[
\varepsilon _{n}=\frac{5}{288n^{3}}-\frac{343}{8640n^{4}}+\frac{2621}{%
41\,472n^{5}}+\frac{2621}{41\,472n^{5}}-\frac{300\,901}{3483\,648n^{6}}.
\]

Being very accurate, we are convinced that the inequalities presented in
Theorem 1 can be succesfully used to obtain other new results.

\begin{acknowledgement}
Computations in this paper were made using Maple software, but they can be
also made (or verified) by direct approach.
\end{acknowledgement}

\begin{acknowledgement}
The work of the second author was supported by a grant of the Romanian
National Authority for Scientific Research, CNCS-UEFISCDI project number
PN-II-ID-PCE-2011-3-0087.
\end{acknowledgement}

\end{document}